\documentclass[a4paper]{article}

\usepackage[utf8]{inputenc}
\usepackage{lmodern}
\usepackage{array}
\usepackage{graphicx}
\usepackage{subfig}
\usepackage{amssymb, amsfonts, amsthm, amsmath}
\usepackage{enumerate}
\usepackage[english]{babel}
\usepackage{tikz}
\usetikzlibrary{decorations}
\usetikzlibrary{decorations.pathreplacing}
\tikzset{individu/.style={draw,thick}}

\numberwithin{equation}{section}

\theoremstyle{plain}
\newtheorem{theorem}{Theorem}[section]
\newtheorem{corollary}[theorem]{Corollary}

\newtheorem{lemma}[theorem]{Lemma}

\theoremstyle{definition}

\theoremstyle{remark}

\newtheorem{example}[theorem]{Example}

\newcommand{\Z}{\mathbb{Z}}

\renewcommand{\epsilon}{\varepsilon}
\renewcommand{\phi}{\varphi}

\newcommand{\Addresses}{{
  \bigskip
  \footnotesize

  Ksenia Chernysh, \textsc{School of Mathematical and Computer Sciences, Heriot-Watt University, Edinburgh, EH14, 4AS, UK}\par\nopagebreak
  \textit{E-mail address}: \texttt{xeniachernysh@gmail.com}
  
  \medskip

  Sanjay Ramassamy, \textsc{Mathematics Department, Brown University, Box 1917, 151 Thayer street, Providence, RI 02912, USA}\par\nopagebreak
  \textit{E-mail address}: \texttt{sanjay\_ramassamy@brown.edu}

}}

\title{Coupling any number of balls in the infinite-bin model}
\author{Ksenia Chernysh \and Sanjay Ramassamy}
\date{\today}

\usepackage{tikz}

\makeatletter \newcommand \listoftodos{\section*{Todo list} \@starttoc{tdo}}
\newcommand\l@todo[2]
{\par\noindent \textit{#2}, \parbox{10cm}{#1}\par} \makeatother

\begin{document}

\maketitle

\begin{abstract}
The infinite-bin model, introduced by Foss and Konstantopoulos in~\cite{FK}, describes the Markovian evolution of configurations of balls placed inside bins, obeying certain transition rules. We prove that we can couple the behaviour of any finite number of balls, provided at least two different transition rules are allowed. This coupling makes it possible to define the regeneration events needed by Foss and Zachary in~\cite{FZ} to prove convergence results for the distribution of the balls.
\end{abstract}

\section{Introduction}
\label{sec:introduction}

\subsection{Background}
\label{subsec:background}

Interacting particle systems and their stochastic evolution have been widely studied in probability and statistical physics. One such system, the infinite-bin model, was introduced by Foss and Konstantopoulos in~\cite{FK} as an abstraction of stochastic ordered graphs, which are a directed version of Erd\"os-R\'enyi graphs that have applications in queuing theory, mathematical ecology and performance evaluation of computer systems (see~\cite{FK} and references therein). A very similar (though less general) model had also appeared in earlier work of Aldous and Pitman in~\cite{AP}.

In the infinite-bin model, a configuration is made up of infinitely many bins indexed by the nonpositive integers, each bin containing a positive and finite number of balls. An elementary random move consists in picking one ball at random (according to a certain probability distribution) and adding one ball to the bin situated immediately to the right of the ball we picked. If the ball we picked was already in the rightmost bin, we create a new bin to its right, add a ball in it and relabel the bins so that the new rightmost bin has label $0$. The stochastic dynamics arises from the iteration of i.i.d. elementary random moves. 

Questions of interest include the existence and uniqueness of a stationary solution, the convergence to the stationary solution for an arbitrary initial configuration and the rate of creation of new bins. The first two questions have been tackled in certain cases by~\cite{FK} and~\cite{FZ}. The last one was discussed in some cases by~\cite{FK} and will be addressed further in an upcoming joint paper with B.~Mallein~\cite{CMR}.

In this paper, we prove that in all the nontrivial cases, we can find a sequence of moves such that, after applying that sequence of moves, the position of the rightmost $N$ balls is some prescribed configuration, independent of the initial configuration (here $N$ is an arbitrary positive integer).

This is a result needed by Foss and Zachary in~\cite{FZ} to prove the convergence to the stationary distribution. They do not prove it in their paper, but instead point to the present paper for a proof. Note that the present paper relaxes the assumptions of~\cite{FZ} so that we can now cover all the cases where the theory of renovation events works.

We finally mention that we can rephrase this result in the framework of automata theory. The infinite-bin model (or rather the finite projections of it defined in the next subsection) can be seen as a finite automaton, and the sequence of moves we construct corresponds to a synchronizing word in the language of automata theory (see e.g.~\cite{Tr}).

\subsection{Definitions and main result}
\label{subsec:definitions}

An \emph{infinite configuration} $X$ is a sequence of positive integers $\left(X(i)\right)_{i\in\Z_-}$ indexed by $\Z_-=\left\{0,-1,-2,\ldots\right\}$. $X(i)$ represents the number of balls in the bin labeled by $i$. We adopt the unusual indexing by nonpositive integers to conform to the original definition of the infinite-bin model in \cite{FK}. It finds its roots in an application to stochastic ordered graphs.

For any integer $k\geq1$, we are going to construct the \emph{move of type $k$} as a function $\phi_k$ from the set of infinite configurations to itself. Fix an infinite configuration $X$. The infinite configuration $Y=\phi_k(X)$ is defined as follows. If $X(0)\geq k$, we set 

\begin{equation}
Y(i)=
\begin{cases}
1 &\mbox{if } i=0\\ 
X(i+1) &\mbox{if } i<0.
\end{cases} 
\end{equation}

If $X(0)<k$, we define

\begin{equation}
i_k=\min\left\{j\in\Z_-\,\middle|\,\sum\limits_{i=j}^0{X(i)}<k\right\}
\end{equation}
and set 

\begin{equation}
Y(i)=
\begin{cases}
X(i)+1 &\mbox{if } i=i_k\\ 
X(i) &\mbox{if } i\neq i_k.
\end{cases} 
\end{equation}

\begin{figure}[htbp]
\centering
\subfloat[A configuration $X$]{\includegraphics[height=1.5in]{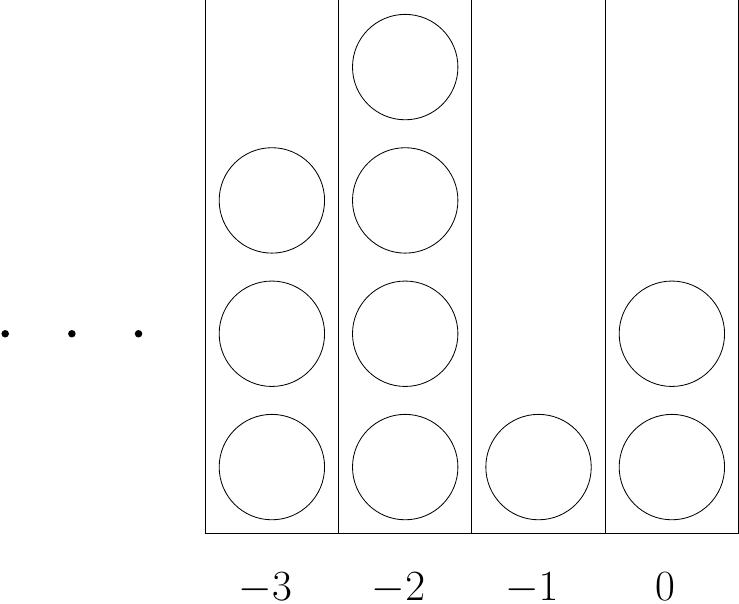}}
\hspace{60pt}
\subfloat[The configuration $\phi_5(X)$]{\includegraphics[height=1.5in]{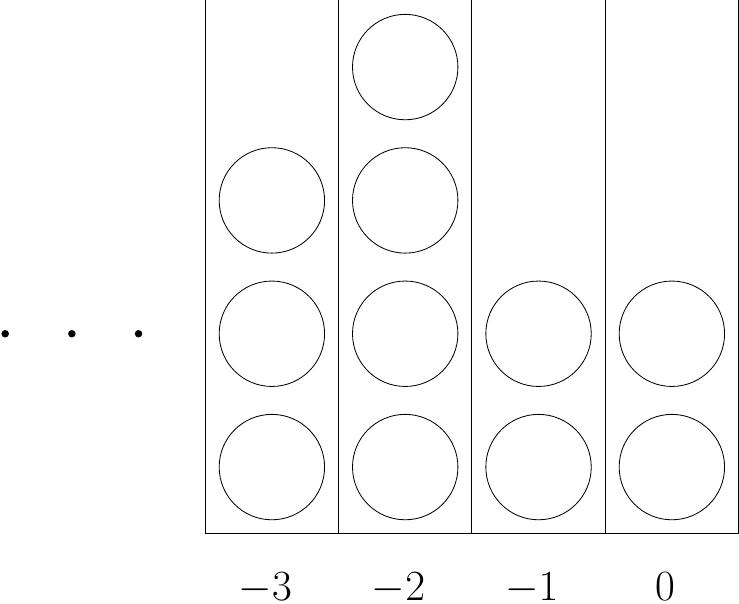}}
\hspace{.2pt}
\subfloat[The configuration $\phi_2(X)$]{\includegraphics[height=1.5in]{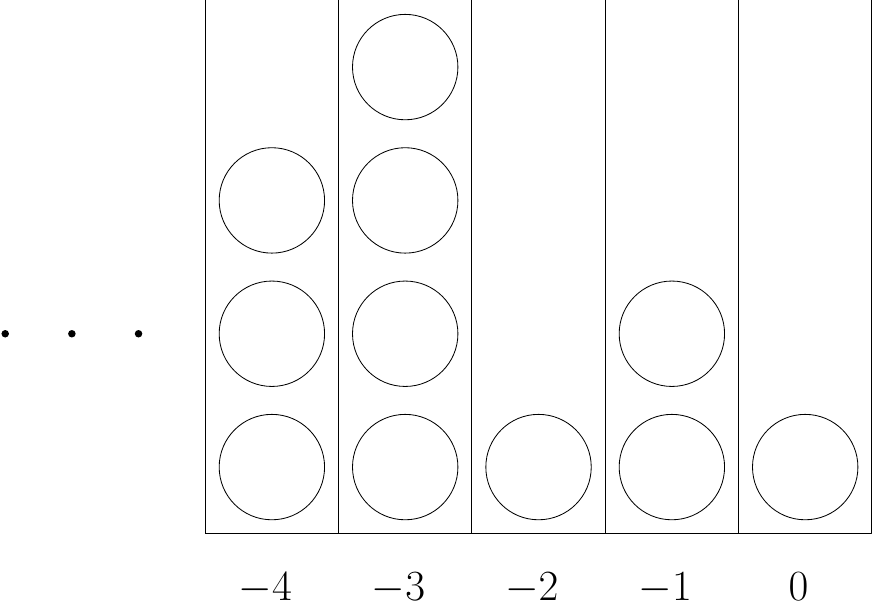}}
\caption{Action of two moves on an infinite configuration.}
\label{fig:move}
\end{figure}

In words, the move of type $k$ adds one ball inside the bin situated immediately to the right of the bin containing the $k$-th ball, where the balls are counted from right to left. When the $k$-th ball is already in the rightmost bin, we create a new bin immediately to its right, add a ball in it and relabel all the bins in such a way that the newly created bin will be labeled by $0$. See Figure~\ref{fig:move} for an example. 

It will be useful to consider configurations with finitely many balls. If $n$ is a positive integer, we define an \emph{$n$-configuration} $X$ to be a sequence of nonnegative integers $\left(X(i)\right)_{i\in\Z_-}$ such that:

\begin{enumerate}
\item $\sum\limits_{i\in\Z_-}{X(i)}=n$
\item there exists an integer $p\geq1$ such that $X(i)>0$ if $i>-p$ and $X(i)=0$ otherwise.
\end{enumerate}

We will then denote the $n$-configuration by the $p$-tuple $\left[X(-p+1),\ldots,X(0)\right]$, omitting the infinite string of $0$'s on the left.

If $X$ is an infinite configuration and $n$ is a positive integer, we define $\pi_n(X)$ the \emph{$n$-ball projection} of $X$ to be the $n$-configuration obtained from $X$ by keeping the rightmost $n$ balls and erasing all the other balls. More precisely,

\begin{equation}
\left(\pi_n(X)\right)_i=
\begin{cases}
\min\left(X(i),n-\sum\limits_{j=i+1}^0{X(j)}\right) &\mbox{if } \sum\limits_{j=i+1}^0{X(j)}<n\\ 
0 &\mbox{otherwise}
\end{cases} 
\end{equation}

An \emph{algorithm} $\Phi$ is defined as the composition of a finite number of moves. It maps one infinite configuration to another. The \emph{length} of an algorithm is the number of moves used to build it. For example, $\Phi=\phi_5^3\phi_2$ is an algorithm of length $4$ obtained by applying $\phi_2$ first followed by three times $\phi_5$.

An algorithm $\Phi$ is said to \emph{couple the first $n$ balls} if for any two infinite configurations $X$ and $Y$, we have $\pi_n(\Phi(X))=\pi_n(\Phi(Y))$. We have the following result regarding the existence of coupling algorithms :

\begin{theorem}
\label{thm:coupling}
For any integers $N\geq1$ and $1\leq k<l$, we can find an algorithm using only the moves $\phi_k$ and $\phi_l$ which couples the first $N$ balls. The length of the algorithm can be taken to be less than $N+4l^2$.
\end{theorem}

The proof is provided in Section~\ref{sec:prooftheorem}. We will see in the proof that we can impose that any infinite configuration obtained after applying that coupling algorithm contains at least $k$ balls in the rightmost bin, which is exactly the setting in which Lemma 2 of \cite{FZ} was stated. Note that we don't need the condition (required by that Lemma 2) that $k$ and $l$ are co-prime.

\subsection{Probabilistic implications}
\label{subsec:probability}

Coupling algorithms are interesting from a probabilistic point of view because they make it possible to define regeneration events, and thus make the infinite-bin model more tractable.

Let $\xi$ be a random variable taking values in the set of the positive integers, and let $(\xi_n)_{n\geq0}$ be a sequence of i.i.d. random variables distributed like $\xi$. Then any sequence of infinite configurations $(X_n)_{n\geq0}$ such that for all $n\geq0$,

\begin{equation}
X_{n+1}=\phi_{\xi_n}\left(X_n\right)
\end{equation}

is a Markov chain.

Theorem 4 in \cite{FZ} proves, under some assumptions, the convergence of the distribution of the number of balls in the rightmost $k$ bins for any positive integer $k$. It is based on their Lemma 2, which is a more restrictive version of our Theorem~\ref{thm:coupling}. For a proof of that lemma, they point to the present paper. Combining our Theorem~\ref{thm:coupling} with Theorem 4 and Proposition 2 of~\cite{FZ}, we obtain:

\begin{corollary}[\cite{FZ}]
\label{thm:convergence}
Assume that $\xi$ has finite expectation and is not constant a.s.. Then for any integer $k\geq0$ and any initial infinite configuration $X_0$, $\left(X_n(-k),\ldots,X_n(0)\right)$ converges to a proper limiting random vector in the total variation norm. Therefore $X_n$ weakly converges to its proper limit.
\end{corollary}

Here ``proper" means the limiting random variable is finite a.s.. Note that if $\xi$ is a.s. constant equal to $c$ then the dynamics is deterministic and ultimately $c$-periodic.

The proof is the same as that of Theorem 4 in \cite{FZ}, replacing their Lemma 2 by our Theorem~\ref{thm:coupling} and constructing regeneration events. A regeneration event $E$ is an event such that at the time when $E$ starts, conditionally on $E$, the future does not depend on what happened before $E$ (here ``depend'' is used in the algebraic sense, not in the probabilistic sense). They define a regeneration event to be a coupling algorithm followed by an infinite sequence of moves with the property that the balls involved in that infinite sequence of moves are only those balls that were created by the coupling algorithm. Such an event will eventually occur a.s.. This implies that if we run the Markov chain with two different initial configurations but with the same sequence $(\xi_n)_{n\geq0}$, the content of the first $k$ bins of the Markov chains will eventually be the same a.s..

An interesting next step would be to find the coupling algorithms that occur most frequently, i.e. for which the ratio ``probability of algorithm'' divided by ``length of algorithm'' is maximal. When the distribution of $\xi$ is the uniform measure on a finite set of integers, this boils down to finding the shortest coupling algorithms, which is a classical question in automata theory (e.g.~\cite{Tr}).

\section{Proof of Theorem~\ref{thm:coupling}}
\label{sec:prooftheorem}

The structure of the proof is as follows. Firstly, we reduce the proof to finding an algorithm which couples the first $k$ balls when applied to $l$-configurations. For this purpose, we will need to define how moves act on $l$-configurations. Secondly, we define a set of $k$ $l$-configurations $\left\{X_0,\ldots,X_{k-1}\right\}$ and two algorithms $\Psi_1$ and $\Psi_2$ using only the moves $\phi_k$ and $\phi_l$ verifying the following properties:

\begin{enumerate}
\item For any $l$-configuration $X$, $\Psi_1(X)\in\left\{X_0,\dots,X_{k-1}\right\}$.
\item For any $0\leq i\leq k-1$, $\Psi_2(X_i)=X_0$.
\end{enumerate}

It will follow that the algorithm $\Psi_2\Psi_1$ couples the first $l$ balls (and \emph{a fortiori} the first $k$ balls) when applied to $l$-configurations.

The case $k=1$ is easy: applying move $\phi_1$ $N$ times couples the first $N$ balls. From now on we will assume that $k\geq2$.

If $n\geq k$ and X,Y are two infinite configurations such that $\pi_n(X)=\pi_n(Y)$, observe that $\pi_{n+1}(\phi_k(X))=\pi_{n+1}(\phi_k(Y))$. So without loss of generality, we may assume that $N=k$: if we find an algorithm coupling the first $k$ balls, applying $N-k$ more times the move $\phi_k$ will couple the first $N$ balls.

We now define projections and moves applied to $n$-configurations. If $n>m$ one defines $\pi_m(X)$ the $m$-ball projection of an $n$-configuration $X$ in the same way as for infinite configurations. The map $\pi_m$ is a function from $n$-configurations into $m$-configurations. If $X$ is an $n$-configuration and $k\leq n$, we define $\phi_k(X)$ to be the following $n$-configuration: starting from $X$, add one ball to the bin situated immediately to the right of the bin containing the $k$-th ball (where the balls are counted from right to left), then delete one ball from the leftmost bin. When the $k$-th ball is already in the rightmost bin, we create a new bin immediately to its right, add a ball in it and relabel all the bins in such a way that the newly created bin will be labeled by $0$. The map $\phi_k$ is a function from $n$-configurations into $n$-configurations. This definition of $\phi_k$ for $n$-configurations parallels the one for infinite configurations, so that moves commute with projections, i.e. if $k\leq n$ and if $X$ is an infinite configuration or an $m$-configuration with $m\geq n$, then

\begin{equation}
\pi_n(\phi_k(X))=\phi_k(\pi_n(X))
\end{equation}

\begin{figure}[htbp]
\centering
\includegraphics[height=3.5in]{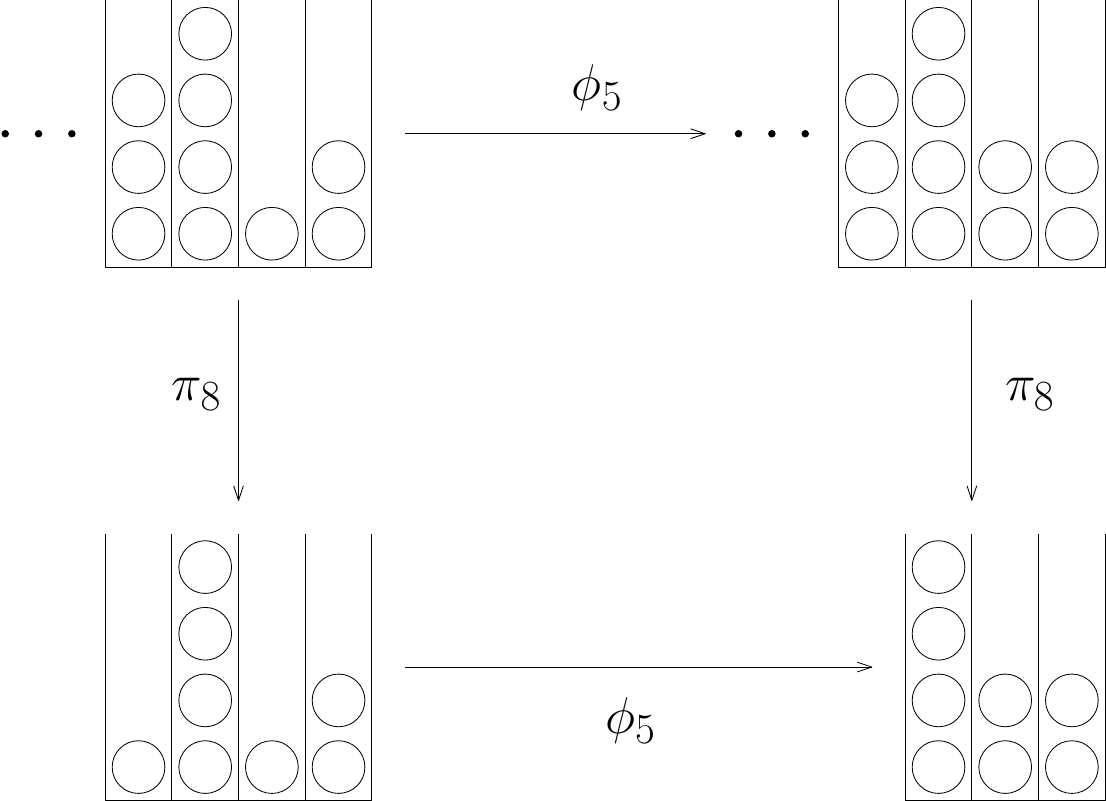}
\caption{Commutation of $\phi_5$ with $\pi_8$.}
\label{fig:commutation}
\end{figure}

This is illustrated in Figure~\ref{fig:commutation}.

Fix an integer $n\geq1$. If, for any move $\phi_k$ entering in the construction of an algorithm $\Phi$, we have $k\leq n$, then $\Phi$ also maps $n$-configurations to $n$-configurations.

Let $\Phi$ be an algorithm using only the moves $\phi_k$ and $\phi_l$ which couples the first $k$ balls when restricted to $l$-configurations. Then if $X$ and $Y$ are two infinite configurations,

\begin{align*}
\pi_k\left(\Phi(X)\right)&=\pi_k\left(\pi_l\left(\Phi(X)\right)\right) \\
&=\pi_k\left(\Phi\left(\pi_l(X)\right)\right) \\
&=\pi_k\left(\Phi\left(\pi_l(Y)\right)\right) \\
&=\pi_k\left(\pi_l\left(\Phi(Y)\right)\right) \\
&=\pi_k\left(\Phi(Y)\right),
\end{align*}

where the third equality follows from the fact that $\Phi$ couples the first $k$ balls of the $l$-configurations $\pi_l(X)$ and $\pi_l(Y)$.

So $\Phi$ also couples the first $k$ balls for infinite configurations. Thus, we can reduce to question to finding an algorithm $\Phi$ which couples the first $k$ balls when applied to $l$-configurations. From now on, by ``configurations" we will mean ``$l$-configurations".

Write $l=kd+r$, where $d$ and $r$ are two integers such that $d\geq1$ and $1\leq r\leq k$.

Define the following configurations. Set $X_0=\left[k,\ldots,k,r\right]$ with $d$ consecutive bins containing $k$ balls and the rightmost bin containing $r$ balls. If $1\leq i\leq r-1$, set $X_i=\left[i,k,\ldots,k,r-i\right]$, with $d$ consecutive bins containing $k$ balls sandwiched between a bin with $i$ balls and a bin with $r-i$ balls. If $r\leq i\leq k-1$, set $X_i=\left[i,k,\ldots,k,k+r-i\right]$, with $d-1$ consecutive bins containing $k$ balls sandwiched between a bin with $i$ balls and a bin with $r-i$ balls. Then for any $1\leq i\leq k$, $\phi_k(X_i)=X_{i-1}$, where indices are taken modulo $k$.

\begin{example}
\label{ex:Xis}
In the case when $k=3$ and $l=11$, we have $d=3$ and $r=2$. Thus, $X_0=\left[3,3,3,2\right]$, $X_1=\left[1,3,3,3,1\right]$ and $X_2=\left[2,3,3,3\right]$.
\end{example}

The rest of the proof of the theorem is based on three lemmas, which will be proved in Section~\ref{sec:prooflemmas}. Set

\begin{equation}
M=\max\left(l,\frac{k(k-1)}{2}\right)
\end{equation}

and $\Psi_1=\phi_k^M$, where $\phi_k^M$ denotes the algorithm built from the move $k$ composed $M$ times. Then we have the following result:

\begin{lemma}
\label{lem:kcycle}
For any $l$-configuration $X$, $\Psi_1(X)\in\left\{X_0,\ldots,X_{k-1}\right\}$.
\end{lemma}

Define the algorithm

\begin{equation}
\Psi=\phi_k^{l-k}\phi_l^{dr+k\frac{d(d-1)}{2}}\phi_k^{k-r}.
\end{equation}

If $n\in\Z$ is an integer, we write $(n)_+=\max(n,0)$. Define the function $f$ from the set $\{0,\ldots,k-1\}$ to itself by:

\begin{equation}
f(i)=\begin{cases}
(k-l+di)_+ &\mbox{if } 0\leq i\leq k-r-1 \\
\left(k-(d+1)(k-i)\right)_+ &\mbox{if } k-r\leq i\leq k-1.
\end{cases}
\end{equation}

Using the function $f$, we can describe the action of $\Psi$ on the configurations $X_i$:

\begin{lemma}
\label{lem:psiaction}
If $0\leq i\leq k-1$, then

\begin{equation}
\Psi(X_i)=X_{f(i)}.
\end{equation}

\end{lemma}

\begin{lemma}
\label{lem:fproperties}
The function $f$ has the following properties:
\begin{enumerate}
\item $f(0)=0$ ;
\item for any $1\leq i\leq k-1$, we have $f(i)<i$.
\end{enumerate}
\end{lemma}

From Lemma~\ref{lem:fproperties}, we deduce that for any $0\leq i\leq k-1$, we have $f^{k-1}(i)=0$. So by Lemma~\ref{lem:psiaction}, for any $0\leq i\leq k$, we have $\Psi^{k-1}(X_i)=X_0$. Combining this with Lemma~\ref{lem:kcycle}, we conclude that for any $l$-configuration $X$, we have $\Psi_2\Psi_1(X)=X_0$, provided we define $\Psi_2=\Psi^{k-1}$.

The configuration $X_0$ has $k$ balls in the rightmost bin, confirming the claim made after the statement of Theorem~\ref{thm:coupling} in Section~\ref{sec:introduction}. Also, note that we did better than coupling the first $k$ balls, we coupled the first $l$ balls.

\begin{example}
Continuing with Example~\ref{ex:Xis}, when $k=3$ and $l=11$ we have $M=11$ and $\Psi=\phi_3^8 \phi_{11}^{15} \phi_3$, thus $\Psi_1=\phi_3^{11}$ and $\Psi_2=\Psi^2=\phi_3^8 \phi_{11}^{15} \phi_3^9 \phi_{11}^{15} \phi_3$. Consider the $11$-configuration $X=\left[2,3,1,3,2\right]$. We first have $\Psi_1(X)=\left[1,3,3,3,1\right]=X_1$, then $\Psi(X_1)=\left[3,3,3,2\right]=X_0$ and $\Psi(X_0)=X_0$, thus $\Psi_2\Psi_1(X)=X_0$.
\end{example}

The length of $\Psi_2\Psi_1$ is

\begin{equation}
L=\max\left(l,\frac{k(k-1)}{2}\right)+\left(k-1\right)\left(k-r+dr+k\frac{d(d+1)}{2}+l-k\right).
\end{equation}

Using that $r\leq k\leq l$ and $d\leq l/k$, we obtain

\begin{align}
L&\leq l+\frac{k(k-1)}{2}+k\left(dr+dk+k\frac{d(d-1)}{2}+l\right) \\
&\leq l+\frac{k^2}{2}+k\left(l+l+k\frac{l^2}{2k^2}+l\right) \\
&\leq l+\frac{l^2}{2}+3kl+\frac{l^2}{2} \\
&\leq l+4l^2.
\end{align}

If we want to couple $N$ balls with $N>l$ for infinite configurations, we will need at most $N-l$ more iterations of $\phi_l$ following $\Phi$. Hence an upper bound of the total length of such a coupling algorithm is $N+4l^2$.

This concludes the proof of the theorem.

\section{Proof of the lemmas}
\label{sec:prooflemmas}

\subsection{Proof of Lemma~\ref{lem:kcycle}}
\label{subsec:prooflemma1}

We will use the following characterization of the configurations $X_i$ : if an $l$-configuration $X$ has at most $k$ balls in each bin, and every bin except maybe for the leftmost and the rightmost bins contains exactly $k$ balls, then $X\in\left\{X_0,\ldots,X_{k-1}\right\}$. Also, in this subsection, it might be convenient to think that the bins don't get relabeled after a new rightmost bin gets created by a move.

Fix an $l$-configuration $X$ and write $\pi_k(X)=\left[i_1,i_2,\ldots,i_p\right]$ for some positive integers $i_1,\ldots,i_p$. After applying $i_1$ times $\phi_k$ to $X$, the bin of $X$ that initially contained $i_2$ balls now contains $i_1+i_2$ balls. Iterating this process, we observe that if we set

\begin{equation}
n=\sum\limits_{j=1}^{p-1}{(p-j)i_j},
\end{equation}

then after applying $n$ times $\phi_k$ to $X$, the rightmost bin of $X$ (which initially contained $i_p$ balls) contains $i_1+\cdots+i_p=k$ balls. All further applications of $\phi_k$ will add new bins to the right of the original rightmost bin of $X$. These new bins will all contain $k$ balls, except maybe for the last (rightmost) bin, which will contain at most $k$ balls.

Note that we can rewrite

\begin{equation}
n=\sum\limits_{i=1}^{k}{(\mbox{distance of the bin containing ball number } i \mbox{ to the rightmost bin})},
\end{equation}

where the balls are counted from right to left. The distance between two bins is given by the absolute value of the difference of their labels. That sum is clearly maximal when $p=k$ and $i_1=\cdots=i_k=1$, thus

\begin{equation}
n\leq\frac{k(k-1)}{2}.
\end{equation}

Recall that we had

\begin{equation}
M=\max\left(l,\frac{k(k-1)}{2}\right).
\end{equation}

Since $M\geq l$, all the balls that were originally in $X$ have been deleted by applying $\phi_k^M$. Since $M\geq n$, after applying $\phi_k^M$, all the bins to the right of the original rightmost bin will contain $k$ balls, except maybe for the leftmost and the rightmost, which will nevertheless contain at most $k$ balls. So $\phi_k^M(X)\in\left\{X_0,\ldots,X_{k-1}\right\}$.

\subsection{Proof of Lemma~\ref{lem:psiaction}}
\label{subsec:prooflemma2}

If $1\leq j\leq l-1$, set $Y_j=\left[j,l-j\right]$ and set $Y_0=\left[l\right]$. For any $1\leq j\leq l$, $\phi_l(Y_j)=Y_{j-1}$, where indices are taken modulo $l$.

We will be using the following lemma, describing how powers of $\phi_k$ and $\phi_l$ act on certain configurations.

\begin{lemma}
\label{lem:klaction}
\begin{enumerate}
\item Assume $i_1,\ldots,i_p$ are $p$ positive integers summing to $l$ and $n=\sum\limits_{j=1}^{p-2}{(p-1-j)i_j}$. Then
\begin{equation}
\phi_l^n\left(\left[i_1,i_2,\ldots,i_p\right]\right)=Y_{l-i_p}.
\end{equation}
\item For any $0\leq j\leq l-1$, we have
\begin{equation}
\phi_k^{l-k}(Y_j)=X_{(k-l+j)_+}.
\end{equation}
\end{enumerate}
\end{lemma}

\emph{Proof of Lemma~\ref{lem:klaction}}. The first part is obtained by iterating the following procedure to remove the leftmost bin:

\begin{equation}
\phi_l^{i_1}\left(\left[i_1,i_2,\ldots,i_p\right]\right)=\left[i_1+i_2,i_3,\ldots,i_p\right].
\end{equation}

For the second part, $Y_j=\left[j,l-j\right]$, with the leftmost bin being empty if $j=0$. We need to distinguish three cases.

\begin{itemize}

\item Assume that $k-l+j\leq0$. Then the rightmost bin already contains at least $k$ balls, so $\phi_k^{l-k}(\left[j,l-j\right])=X_0=X_{(k-l+j)_+}$.

\item Assume that $0<k-l+j\leq l-k$. Then $k-l+j$ and $2l-2k-j$ are both nonnegative. We first compute

\begin{equation}
\phi_k^{k-l+j}(\left[j,l-j\right])=\left[l-k,k\right].
\end{equation}

Since $2l-2k-j\geq l-2k$, applying $\phi_k^{2l-2k-j}$ to $\left[l-k,k\right]$ will delete at least $l-2k$ balls from the leftmost bin (initially containing $l-k$ balls) and will construct columns of size $k$, except the rightmost column which will contain at most $k$ balls. So $\phi_k^{2l-2k-j}(\left[l-k,k\right])$ has to be equal to some $X_i$. The number of balls in the rightmost bin of $X_i$ will be congruent to $2l-2k-j=(2d-2)k+2r-j$ modulo $k$, so the number of balls in its leftmost bin will be congruent to $j-r=k-l+j+(d-1)k$ modulo $k$. So $\phi_k^{2l-2k-j}(\left[l-k,k\right])=X_{k-l+j}$. So

\begin{equation}
\phi_k^{l-k}(\left[j,l-j\right])=X_{k-l+j}=X_{(k-l+j)_+}.
\end{equation}

\item Assume that $k-l+j>l-k$. Then

\begin{equation}
\phi_k^{l-k}(\left[j,l-j\right])=\left[j+k-l,2l-k-j\right],
\end{equation}
because the rightmost bin will always have less than $k$ balls in the process of applying $l-k$ times $\phi_k$. Moreover, the leftmost bin of $\left[j+k-l,2l-k-j\right]$ also has less than $k$ balls, so that configuration has to be one of the $X_i$'s, namely $X_{k-l+j}$. Thus, $\phi_k^{l-k}(\left[j,l-j\right])=X_{(k-l+j)_+}$.

\end{itemize}

This concludes the proof of Lemma~\ref{lem:klaction}.

Let us now use it to prove Lemma~\ref{lem:psiaction}. We need to distinguish two cases.

\begin{itemize}

\item Assume that $0\leq i\leq k-r-1$. Then $\phi_k^{k-r}(X_i)=X_{i+r}$, with $r\leq i+r\leq k-1$. So $X_{i+r}=\left[i+r,k,\ldots,k,k-i\right]$, with $d-1$ bins containing $k$ balls sandwiched between a bin with $i+r$ balls and a bin with $k-i$ balls. If we set

\begin{equation}
m_1=(i+r)(d-1)+\sum\limits_{j=1}^{d-2}{kj},
\end{equation}

then by the first part of Lemma~\ref{lem:klaction} we get

\begin{equation}
\phi_l^{m_1}\left(\left[i+r,k,\ldots,k,k-i\right]\right)=Y_{l-k+i}.
\end{equation}

Set 

\begin{equation}
m_2=dr+k\frac{d(d-1)}{2}-m_1=l-k+i-di.
\end{equation}

Then we get

\begin{equation}
\phi_l^{m_2}\left(Y_{l-k+i}\right)=Y_{di}.
\end{equation}

By the second part of Lemma~\ref{lem:klaction}, we obtain

\begin{equation}
\phi_k^{l-k}(Y_{di})=X_{(k-l+di)_+}.
\end{equation}

Putting everything together, we conclude that $\Psi(X_i)=X_{f(i)}$ in this case.

\item Assume that $k-r\leq i\leq k-1$. Then $\phi_k^{k-r}(X_i)=X_{i+r-k}$ with $0\leq i+r-k\leq r-1$. So $X_{i+r-k}=\left[i+r-k,k,\ldots,k,k-i\right]$ with $d$ bins containing $k$ balls sandwiched between a bin with $i+r-k$ balls and a bin with $k-i$ balls. If we set

\begin{equation}
n_1=d(i+r-k)+\sum\limits_{j=1}^{d-1}{kj},
\end{equation}

then by the first part of Lemma~\ref{lem:klaction} we get

\begin{equation}
\phi_l^{n_1}(\left[i+r-k,k,\ldots,k,k-i\right])=Y_{l-k+i}.
\end{equation}

Set

\begin{equation}
n_2=dr+k\frac{d(d-1)}{2}-n_1=d(k-i).
\end{equation}

Then we get

\begin{equation}
\phi_l^{n_2}\left(Y_{l-k+i}\right)=Y_{r-k+i(d+1)}.
\end{equation}

By the second part of Lemma~\ref{lem:klaction}, we obtain

\begin{equation}
\phi_k^{l-k}(Y_{r-k+i(d+1)})=X_{(k-l+r-k+i(d+1))_+}=X_{\left(k-(d+1)(k-i)\right)_+}.
\end{equation}

Putting everything together, we conclude that $\Psi(X_i)=X_{f(i)}$ in this case too.

\end{itemize}

\subsection{Proof of Lemma~\ref{lem:fproperties}}
\label{subsec:prooflemma3}

Firstly, $f(0)=(k-l)_+=0$ because $k<l$. Now let us show by descending induction on $i$ that for $1\leq i\leq k-1$, we have $f(i)<i$. Using that $d\geq1$ and $k\geq2$, we have $f(k-1)=(k-d-1)_+ < k-1$. Fix $1\leq i\leq k-2$. If $f(i)=0$, we are done, so let us assume that $f(i)>0$. It is easy to check that $f$ is nondecreasing, thus $f(i+1)>0$. Considering the three cases where $i$ is less than, equal to or greater than $k-r-1$, one proves that $f(i+1)-f(i)\geq d\geq1$. Using the induction hypothesis, we obtain that $f(i)\leq f(i+1)-1< i+1-1$, which concludes the proof of the second part of the lemma.

\section*{Acknowledgements}

We thank Sergey Foss for fruitful discussions. SR acknowledges the hospitality of the Mathematical Sciences Research Institute in Berkeley during the ``Random spatial processes" program, where this work was started.

\label{Bibliography}
\bibliographystyle{plain}
\bibliography{bibliographie}

\Addresses

\end{document}